# HOW LARGE A DISC IS COVERED BY A RANDOM WALK IN $N$ STEPS?

BY AMIR DEMBO,[1] YUVAL PERES[2] AND JAY ROSEN[3]

*Stanford University, University of California at Berkeley and College of Staten Island, CUNY*

We show that the largest disc covered by a simple random walk (SRW) on $\mathbb{Z}^2$ after $n$ steps has radius $n^{1/4+o(1)}$, thus resolving an open problem of Révész [*Random Walk in Random and Non-Random Environments* (1990) World Scientific, Teaneck, NJ]. For any fixed $\ell$, the largest disc completely covered at least $\ell$ times by the SRW also has radius $n^{1/4+o(1)}$. However, the largest disc completely covered by each of $\ell$ independent simple random walks on $\mathbb{Z}^2$ after $n$ steps is only of radius $n^{1/(2+2\sqrt{\ell})+o(1)}$. We complement this by showing that the radius of the largest disc completely covered at least a fixed fraction $\alpha$ of the maximum number of visits to any site during the first $n$ steps of the SRW on $\mathbb{Z}^2$, is $n^{(1-\sqrt{\alpha})/4+o(1)}$. We also show that almost surely, for infinitely many values of $n$ it takes about $n^{1/2+o(1)}$ steps after step $n$ for the SRW to reach the first previously unvisited site (and the exponent $1/2$ is sharp). This resolves a problem raised by Révész [*Ann. Probab.* **21** (1993) 318–328].

**1. Introduction.** Consider the simple random walk (SRW) on $\mathbb{Z}^2$ starting at the origin and run for $n$ steps. Let $\mathbf{R}_n$ denote the radius of the largest disc centered at the origin that is completely covered by the walk (throughout this paper, "disc" refers to the intersection of $\mathbb{Z}^2$ with a Euclidean disc, but all our results apply if one takes a square instead). In [3], Theorem 1.4 we showed that for all $y > 0$,

$$(1.1) \qquad \lim_{n \to \infty} \mathbf{P}\left(\frac{(\log \mathbf{R}_n)^2}{\log n} \geq y\right) = e^{-4y},$$

Received March 2005; revised April 2006.

[1]Supported in part by NSF Grants DMS-04-06042 and DMS-FRG-02-44323.

[2]Supported in part by NSF Grants DMS-01-04073 and DMS-02-44479.

[3]Supported in part by grants from the NSF and from PSC-CUNY.

*AMS 2000 subject classifications.* Primary 60G50; secondary 60G17, 82C41.

*Key words and phrases.* Planar random walk, favorite points, covered discs.

as conjectured by Kesten and Révész.

If we ask for the largest disc covered after $n$ steps of the SRW without specifying the center of the disc, then the answer changes dramatically.

THEOREM 1.1.   *If $\widetilde{\mathcal{R}}(n)$ denotes the radius of the largest disc completely covered by a SRW on $\mathbb{Z}^2$ after $n$ steps, then almost surely $\widetilde{\mathcal{R}}(n) = n^{1/4+o(1)}$, that is,*

$$\lim_{n \to \infty} \frac{\log \widetilde{\mathcal{R}}(n)}{\log n} = \frac{1}{4} \qquad a.s. \tag{1.2}$$

*Alternatively, with $\mathcal{R}(r)$ denoting the radius of the largest disc completely covered by a SRW on $\mathbb{Z}^2$ before its first exit of $D(0,r) = \{x \in \mathbb{Z}^2 : |x| < r\}$,*

$$\lim_{r \to \infty} \frac{\log \mathcal{R}(r)}{\log r} = \frac{1}{2} \qquad a.s. \tag{1.3}$$

The problem of finding the radius of the largest disc completely covered by a SRW on $\mathbb{Z}^2$ after $n$ steps was first raised by Révész [8], page 247, who later found upper and lower bounds for the ratio $\log \widetilde{\mathcal{R}}(n)/\log n$ (see [9]). We thank Zhan Shi for informing us of simulations by Arvind Singh which indicated that this ratio tends to 1/4.

If we require that our disc be multiply covered we obtain the following.

THEOREM 1.2.   *If $\widetilde{\mathcal{R}}(n;k)$ denotes the radius of the largest disc completely covered at least $k$ times by a SRW on $\mathbb{Z}^2$ after $n$ steps, then for any $0 < \alpha < 1$,*

$$\lim_{n \to \infty} \frac{\log \widetilde{\mathcal{R}}(n; \alpha(\log n)^2/\pi)}{\log n} = \frac{1 - \sqrt{\alpha}}{4} \qquad a.s. \tag{1.4}$$

*Consequently, for any fixed $k \geq 1$,*

$$\lim_{n \to \infty} \frac{\log \widetilde{\mathcal{R}}(n; k)}{\log n} = \frac{1}{4} \qquad a.s. \tag{1.5}$$

*Alternatively, with $\mathcal{R}(r;k)$ denoting the radius of the largest disc completely covered at least $k$ times by a SRW in $\mathbb{Z}^2$ before its first exit of $D(0,r)$, we have that for any $0 < \alpha < 1$*

$$\lim_{r \to \infty} \frac{\log \mathcal{R}(r; 4\alpha(\log r)^2/\pi)}{\log r} = \frac{1 - \sqrt{\alpha}}{2} \qquad a.s. \tag{1.6}$$

We note in passing that (1.4) deals with the largest disc of $\alpha$-favorite sites for the SRW on $\mathbb{Z}^2$ by time $n$, whereas [2], Section 5 provides information about the number of such sites. Further, since $\widetilde{\mathcal{R}}(n) \geq \widetilde{\mathcal{R}}(n;k) \geq$ $\widetilde{\mathcal{R}}(n; \alpha(\log n)^2/\pi)$ for all $n$ sufficiently large, the statement (1.5) is an immediate consequence of (1.4) and (1.2).

We can generalize Theorem 1.1 by considering $\ell$ independent simple random walks on $\mathbb{Z}^2$.

THEOREM 1.3.  *If $\widetilde{\mathcal{R}}_\ell(n)$ denotes the radius of the largest disc completely covered by each of $\ell$ independent simple random walks on $\mathbb{Z}^2$ after $n$ steps, then*

$$\lim_{n \to \infty} \frac{\log \widetilde{\mathcal{R}}_\ell(n)}{\log n} = \frac{1}{2 + 2\sqrt{\ell}} \qquad a.s. \tag{1.7}$$

*Alternatively, with $\mathcal{R}_\ell(r)$ denoting the radius of the largest disc completely covered by each of $\ell$ independent SRWs on $\mathbb{Z}^2$, each of whom is run until it first exits $D(0, r)$,*

$$\lim_{r \to \infty} \frac{\log \mathcal{R}_\ell(r)}{\log r} = \frac{1}{1 + \sqrt{\ell}} \qquad a.s. \tag{1.8}$$

The SRW needs about $r^2$ steps to exit a disc of radius $r$. Thus, considering the random times $n$ in which the SRW is sufficiently inside a completely covered disc of radius roughly $n^{1/4}$, with Theorem 1.1 we can also solve a related problem raised by Révész in [9].

THEOREM 1.4.  *If $V(n)$ is the number of steps after step $n$ until the SRW on $\mathbb{Z}^2$ first visits any of the previously unvisited sites, then*

$$\limsup_{n \to \infty} \frac{\log V(n)}{\log n} = \frac{1}{2} \qquad a.s. \tag{1.9}$$

Of course $\liminf_{n \to \infty} V(n) = 1$.

We note in passing that the situation is quite different for the SRW on $\mathbb{Z}^d$ for $d \geq 3$, where due to the transience of the process, one has that

$$\lim_{n \to \infty} \frac{\log \widetilde{\mathcal{R}}(n)}{\log \log n} = \frac{1}{d - 2} \qquad a.s. \tag{1.10}$$

as shown in [5], and for $d = 1$, where $2\widetilde{\mathcal{R}}(n)$ is the difference between the maximum and the minimum of the SRW, which upon scaling by $n^{-1/2}$ converges in law to an explicit nondegenerate limit (for finer information on the favorite sites in one dimension, cf. [11] and the references therein).

We now explain the intuitive picture behind our results, starting with Theorem 1.1. To this end, let $\tau(r)$ denote the number of steps until the SRW first exits $D(0, r)$. Clearly, $\mathcal{R}(r) = \widetilde{\mathcal{R}}(\tau(r))$ and (1.2) is equivalent to (1.3) since

$$\lim_{r \to \infty} \frac{\log \tau(r)}{\log r} = 2 \qquad a.s. \tag{1.11}$$

Similarly, (1.11) implies the equivalence between (1.4) and (1.6), as well as the equivalence between (1.7) and (1.8).

Turning to (1.3), let $r_k = e^k$ and for any $x \in \mathbb{Z}^2$ consider the family of Euclidean discs centered at $x$, $\{D(x, r_k); k = 0, 1, \ldots, m\}$. Fixing $0 < \gamma < 1$ such that $\gamma m$ is an integer, let $\beta = 1 - \gamma$. Starting from the sphere $\partial D(x, r_{\gamma m})$, the probability of a planar Brownian motion reaching $\partial D(x, r_{\gamma m-1})$ before exiting $D(x, r_m)$ is $1 - 1/(\beta m)$, so that the probability of it making $a\beta^2 m^2$ excursions from $\partial D(x, r_{\gamma m})$ to $\partial D(x, r_{\gamma m-1})$ before exiting $D(x, r_m)$ (referred to as $\gamma$-excursions) is about $e^{-a\beta m}$. Since there are about $e^{2\beta m}$ disjoint discs of radius $r_{\gamma m}$ in $D(0, r_m)$, ignoring the fact that they have different centers, we expect that the maximal number of such $\gamma$-excursions among all discs of radius $r_{\gamma m}$ is about $2\beta^2 m^2$. Further, the probability of the Brownian motion not hitting the disc $D(x, r_0)$ during one $\gamma$-excursion is $1 - 1/(\gamma m)$, so that the probability of not hitting it during $a\gamma^2 m^2$ consecutive $\gamma$-excursions is about $e^{-a\gamma m}$.

Suppose for the moment that the same applies for the SRW, namely, the maximal number of $\gamma$-excursions that the SRW makes is $2\beta^2 m^2$ and the probability of it not hitting the center of a disc during $a\gamma^2 m^2$ consecutive $\gamma$-excursions is about $e^{-a\gamma m}$. Then, since there are about $e^{2\gamma m}$ points in each $D(x, r_{\gamma m})$, ignoring the fact that they are not centered, the expected number of points not visited during $a\gamma^2 m^2$ consecutive $\gamma$-excursions is about $e^{(2-a)\gamma m}$. Hence about $2\gamma^2 m^2$ of the $\gamma$-excursions are needed for the SRW to visit all sites in a given disc $D(x, r_{\gamma m})$. To find the maximal value of $\gamma$ for which some disc of radius $r_{\gamma m}$ is covered by the SRW, equate $2\beta^2 m^2$ and $2\gamma^2 m^2$ (for $\beta = 1 - \gamma$), to get $\gamma = \beta = 1/2$, as stated.

Also, by the preceding reasoning the probability of each of $\ell$ independent random walks having $a\beta^2 m^2$ or more $\gamma$-excursions for a given disc $D(x, r_{\gamma m})$ is $e^{-a\ell\beta m}$. Considering all possible discs, the maximal value of the preceding parameter $a$ is $2/\ell$. As about $2\gamma^2 m^2$ of the $\gamma$-excursions are needed for each SRW to visit all sites in such a disc, one equates $(2/\ell)\beta^2 m^2$ and $2\gamma^2 m^2$ to get in the context of Theorem 1.3 that $\gamma = 1/(1 + \sqrt{\ell})$, as stated.

To predict the result of Theorem 1.2 one uses the same reasoning, except for replacing the probability of about $\exp(-a\gamma m)$ of the SRW not hitting the center of a disc during $a\gamma^2 m^2$ consecutive $\gamma$-excursions with the probability of about $\exp(-(\sqrt{a} - \sqrt{2\alpha}/\gamma)^2 \gamma m)$ that it makes less than $(4\alpha/\pi)m^2$ visits to the center of the disc during $a\gamma^2 m^2$ $\gamma$-excursions (cf. Lemma 5.1 for the argument leading to this tail probability). Indeed, using the latter probability we find that about $2(\gamma + \sqrt{\alpha})^2 m^2$ of the $\gamma$-excursions are needed to assure that all lattice sites in a given disc $D(x, r_{\gamma m})$ are $\alpha$-favorites. Equating this with $2\beta^2 m^2$ yields the value of $\gamma = (1 - \sqrt{\alpha})/2$ as stated.

To prove these three theorems one needs nontrivial modifications of the classical second moment method. Fortunately, adapting the "multiscale refinement" machinery of [2, 3, 4] to the present context, provides the necessary ingredients for proving these results. Indeed, in Sections 2 and 3 we prove the bounds on $\mathcal{R}(r)$ for Theorem 1.1 (lower bounds and upper bounds, resp.), and in Sections 5 and 6, extend these bounds to the setting of Theorem 1.2, whereas in Section 4 we extend both bounds to the setting of Theorem 1.3. Finally, Section 7 is devoted to the derivation of Theorem 1.4.

When actually proving these theorems one uses potential theory estimates for the SRW (e.g., from [7]), in order to justify that its hitting probabilities are sufficiently close to those of the planar Brownian motion. This, as well as the control of the effect of noncentering and the control of the influence of terminal points of each excursion (when considering the relevant second moment), require a better separation of scales. Thus, the proof is carried along the preceding reasoning, but for a sequence $r_k$ which grows at rate $e^{ck \log k}$ for some $c > 0$ large enough ($c = 3$ suffices here). Specifically, taking throughout $r_k = (k!)^3$ allows us to best reuse proofs from [2, 3, 4], focusing in this manuscript on those ingredients that are not already present there. For the same reason, when proving the lower bounds we also consider $r_{m,k} = r_m/r_k$, so, for example, $r_{m,[\beta m]}$ is roughly of same size as $r_{[\gamma m]}$ for $\gamma = 1 - \beta$ (both being about $e^{c\gamma m \log m}$).

Though we deal here exclusively with the SRW on $\mathbb{Z}^2$, similar results apply for the whole class of random walks considered in [2], Theorem 5.1 upon appropriately modifying the relevant proofs.

The results of this manuscript also inspired the analogous treatment of extremal points of the discrete Gaussian free field in the box $[-n, n]^2$ subject to zero boundary conditions, where, for example, the size of the largest sub-box of $\alpha$-high points of the field corresponds to (1.4) here (see [1] for details).

Throughout this paper we use $o(1_m)$ to denote a function $f(m)$ which converges to zero as $m \to \infty$ and use the notation $a_m \sim b_m$ to indicate that $a_m/b_m \to 1$ as $m \to \infty$.

## 2. The lower bound for Theorem 1.1.

Let $(S_i, i \geq 0)$ denote the SRW on $\mathbb{Z}^2$ with $D(x, r) = \{y \in \mathbb{Z}^2 : |y - x| < r\}$ denoting the disc of radius $r$ centered at $x$. For any set $A \subseteq \mathbb{Z}^2$ we let

$$\partial A = \left\{ y \in \mathbb{Z}^2 : y \notin A, \text{ and } \inf_{x \in A} |y - x| = 1 \right\}$$

denote the boundary of $A$ in $\mathbb{Z}^2$ and $T_A = \inf\{i \geq 0 : S_i \in A\}$ the hitting time of $A$, so in particular $\tau(r) = T_{\partial D(0,r)}$. As in the Introduction, we let $r_0 = 1$ and $r_k = (k!)^3$ for $k \geq 1$. Using the monotonicity of $r \mapsto \mathcal{R}(r)$ and the fact that $\lim_{m \to \infty} \log r_m / \log r_{m-1} = 1$, it follows by a simple interpolation argument that (1.3) is an immediate consequence of

$$(2.1) \qquad \lim_{m \to \infty} \frac{\log \mathcal{R}(r_m)}{\log r_m} = \frac{1}{2} \qquad \text{a.s.}$$

We proceed to provide here the relevant lower bound

$$\liminf_{m\to\infty} \frac{\log \mathcal{R}(r_m)}{\log r_m} \geq \frac{1}{2} \qquad \text{a.s.} \tag{2.2}$$

deferring the corresponding upper bound to Section 3. It actually suffices to prove that for any $\eta > 0$ there exists $p_\eta > 0$ such that for all sufficiently large $m$,

$$\mathbf{P}\left(\frac{\log \mathcal{R}(8r_m)}{\log r_m} \geq \frac{1}{2} - 2\eta\right) \geq p_\eta > 0. \tag{2.3}$$

Indeed, then for all $m$ large enough,

$$\mathbf{P}\left(\frac{\log \mathcal{R}(r_m)}{\log r_{m+1}} \leq \frac{1}{2} - 3\eta\right) \leq 1 - p_\eta < 1.$$

Further, with $\widetilde{\mathcal{R}}([s,t])$ denoting the radius of the largest disc completely covered by $\{S_i : i = s, \ldots, t\}$, we have that

$$\mathcal{R}(r_{m+1}) \geq \max\{\widetilde{\mathcal{R}}([\tau((k-1)r_m), \tau(kr_m)]) : k = 1, \ldots, (m+1)^3\}.$$

So, by the strong Markov property of the SRW at the successive stopping times $\tau(kr_m)$, $k = 1, 2, \ldots, (m+1)^3$ and the fact that $D(S_{\tau((k-1)r_m)}, r_m) \subseteq D(0, kr_m) \subseteq D(0, r_{m+1})$, we get that

$$\mathbf{P}\left(\frac{\log \mathcal{R}(r_{m+1})}{\log r_{m+1}} \leq \frac{1}{2} - 3\eta\right) \leq (1 - p_\eta)^{(m+1)^3}.$$

Consequently, an application of the Borel–Cantelli lemma, followed by taking $\eta \downarrow 0$ yields the bound of (2.2).

Turning to the proof of (2.3), we next construct a subset of the event appearing in (2.3), the probability of which is easier to bound below. To this end, let $r_{m,k} = r_m/r_k$ for $k = 1, \ldots, m$ [so that $r_{m,1} = r_m = (m!)^3$ and $r_{m,m} = 1$]. Then, fixing $a > 0$ we set $n_k = n_k(a) = 3ak^2 \log k$ for $3 \leq k \leq m-1$ and for any $x \in \mathbb{Z}^2$, let $N^x_{m,k}$ denote the number of excursions of $\{S_i\}$ from $\partial D(x, r_{m,k})$ to $\partial D(x, r_{m,k-1})$ until time $T_{\partial D(x,r_m)}$. Fixing $0 < \beta < 1$, with some abuse of notation we let $\beta m$ denote hereafter the integer part of $\beta m$. Let $H^x_{\beta m}$ denotes the event that the SRW visits each point in $D(x, r_{m,\beta m+1})$ during the first $N^x_{m,\beta m}$ excursions from $\partial D(x, r_{m,\beta m})$ to $\partial D(x, r_{m,\beta m-1})$.

We say that a point $x \in \mathbb{Z}^2$ is $m, \beta$-*successful* if

$$H^x_{\beta m} \text{ occurs} \quad \text{and} \quad n_k(a) - k \leq N^x_{m,k} \leq n_k(a) + k \qquad \text{for } k = 3, \ldots, \beta m. \tag{2.4}$$

Let $\mathcal{A}_m \subseteq \mathbb{Z}^2$ be a maximal collection of points in $[3r_m, 4r_m]^2$ such that the distance between any two points in $\mathcal{A}_m$ is at least $4r_{m,\beta m}$.

The existence of an $m, \beta$-successful point in $\mathcal{A}_m$ implies that $\mathcal{R}(8r_m) \geq r_{m,\beta m+1}$. Noting that $\log r_{m,\beta m+1} \sim (1-\beta) \log r_m$, we thus establish (2.3) by showing that for all $\beta = 1 - \gamma > 1/2$, upon taking $a = a(\beta) < 2$ such that

$$(2.5) \qquad a\beta^2 > 2\gamma^2,$$

the probability that there is at least one $m, \beta$-successful point in $\mathcal{A}_m$ is bounded away from zero as $m \to \infty$.

As we will show, the fact that $a < 2$ guarantees that, with a probability that is bounded away from zero as $m \to \infty$, there exists at least one $x \in \mathcal{A}_m$ for which $n_k(a) - k \leq N_{m,k}^x \leq n_k(a) + k$, $k = 3, \ldots, \beta m$, and the relation (2.5) guarantees that with very high probability $H_{\beta m}^x$ then holds as well.

Let $V_m = \sum_{x \in \mathcal{A}_m} Y(m,x)$, where $Y(m,x)$ denotes the indicator random variable for the event $\{x \text{ is } m, \beta\text{-successful}\}$. Then, we have (2.3) as soon as we show that for any $\delta > 0$ and all $m$ sufficiently large,

$$(2.6) \qquad \mathbf{P}(V_m \geq r_{\beta m}^{2-a-\delta}) \geq c_\delta > 0$$

Note that $|\mathcal{A}_m| = r_m^2/(16 r_{m,\beta m}^2) = r_{\beta m}^2/16$ so that by (2.9) of Lemma 2.1, for some $\delta'_m \to 0$,

$$(2.7) \qquad \mathbb{E}(V_m) \geq |\mathcal{A}_m| \bar{q}_m \geq r_{\beta m}^{2-a-\delta'_m}.$$

Applying the Paley–Zygmund inequality (see [6], page 8), it thus suffices to show that $\mathbb{E}(V_m^2) \leq C(\mathbb{E}V_m)^2$ for some $C < \infty$ and all $m$ sufficiently large. Furthermore, $\mathbb{E}V_m \to \infty$ when $m \to \infty$ [see (2.7)], so it suffices to show that

$$(2.8) \qquad \mathbb{E}\left(\sum_{\substack{x,y \in \mathcal{A}_m \\ x \neq y}} Y(m,x)Y(m,y)\right) \leq C(\mathbb{E}V_m)^2.$$

The next lemma is proven at the end of this section.

LEMMA 2.1. *If $a$ and $\beta = 1 - \gamma$ satisfy (2.5) then there exists $\delta_m \to 0$ such that*

$$(2.9) \qquad \bar{q}_m := \inf_{x \in \mathcal{A}_m} \mathbf{P}(x \text{ is } m, \beta\text{-successful}) \geq r_{\beta m}^{-(a+\delta_m)}.$$

*Further, for some $c < \infty$, all $m$ and $x \neq y \in \mathcal{A}_m$,*

$$(2.10) \qquad \mathbf{P}(x, y \text{ are } m, \beta\text{-successful}) \leq c\bar{q}_m^2 r_{k(x,y)}^{a+\delta_{k(x,y)}},$$

*where $k(x,y) = \min\{j \geq 1 : D(x, r_{m,j} + 1) \cap D(y, r_{m,j} + 1) = \varnothing\}$ and $k(x,y) \leq \beta m$ when $x \neq y \in \mathcal{A}_m$.*

We return to the proof of (2.3). In the sequel, we let $C_i$ denote finite constants that are independent of $m$. The definition of $k(x,y) \geq 1$ implies that $|x - y| < 2(r_{m,k(x,y)-1} + 1)$. Note that there are at most $C_0 r_{m,k-1}^2/r_{m,\beta m}^2 = C_0 r_{\beta m}^2/r_{k-1}^2 = C_0' |\mathcal{A}_m| r_{k-1}^{-2}$ points $y \in \mathcal{A}_m$ in the ball of radius $2(r_{m,k-1} + 1)$ centered at $x$. Thus, it follows from Lemma 2.1 that for any fixed $\eta > 0$ such that $2 - (a + 2\eta) > 0$,

$$\sum_{\substack{x,y \in \mathcal{A}_m \\ 1 \leq k(x,y) \leq \beta m}} \mathbb{E}(Y(m,x)Y(m,y))$$

$$\leq C_1 \sum_{\substack{x,y \in \mathcal{A}_m \\ 1 \leq k(x,y) \leq \beta m}} \bar{q}_m^2 r_{k(x,y)}^{a+\eta}$$

(2.11)

$$\leq C_2 \bar{q}_m^2 |\mathcal{A}_m|^2 \sum_{k=1}^{\beta m} r_{k-1}^{-2} r_k^{a+\eta} \leq C_3 (|\mathcal{A}_m| \bar{q}_m)^2 \sum_{k=1}^{\beta m} r_k^{-(2-a-2\eta)}$$

$$\leq C_3 (|\mathcal{A}_m| \bar{q}_m)^2 \sum_{k=1}^{\infty} r_k^{-(2-a-2\eta)} \leq C_4 (\mathbb{E}V_m)^2,$$

which completes the proof of (2.3).

PROOF OF LEMMA 2.1. We say that a point $x \in \mathcal{A}_m \subseteq \mathbb{Z}^2$ is $m, \beta$-*pre-successful* if

$$n_k(a) - k \leq N_{m,k}^x \leq n_k(a) + k \qquad \text{for } k = 3, \ldots, \beta m.$$

The proof of Lemma 3.2 of [10] establishes the analog of the statements (2.9) and (2.10) with $m, \beta$-successful replaced by $m, \beta$-presuccessful, $\beta = 1$ and where instead of $r_{m,k}$ we have $m^{3(m-k)} e^m$. This proof works just as well for our choice of $\beta$ and $r_{m,k}$. Since an $m, \beta$-successful point is also $m, \beta$-presuccessful this establishes the upper bound of (2.10). It thus remains only to show that uniformly in $x \in \mathcal{A}_m$

(2.12) $\mathbf{P}(x \text{ is } m, \beta\text{-successful}) \geq (1 + o(1_m)) \mathbf{P}(x \text{ is } m, \beta\text{-presuccessful}).$

To this end, let $\mathcal{L}_{z,m,\beta}$ denote the event that $z$ is not visited during the first $n_{\beta m}(a) - \beta m$ excursions from $\partial D(x, r_{m,\beta m})$ to $\partial D(x, r_{m,\beta m-1})$. We first show that

(2.13) $\quad \xi_m = \sup_{x \in \mathcal{A}_m} \mathbf{P}\left( \bigcup_{z \in D(x, r_{m,\beta m+1})} \mathcal{L}_{z,m,\beta} \right) \to 0 \qquad \text{as } m \to \infty.$

To see this, note that for $x \in \mathcal{A}_m$

(2.14) $\quad \mathbf{P}\left( \bigcup_{z \in D(x, r_{m,\beta m+1})} \mathcal{L}_{z,m,\beta} \right) \leq 4 r_{m,\beta m+1}^2 \sup_{z \in D(x, r_{m,\beta m+1})} \mathbf{P}(\mathcal{L}_{z,m,\beta}),$

and by the strong Markov property of the SRW,

(2.15) $\quad \mathbf{P}(\mathcal{L}_{z,m,\beta}) \leq \left( \sup_{y \in \partial D(x, r_{m,\beta m})} \mathbf{P}^y(T_z > T_{\partial D(x, r_{m,\beta m-1})}) \right)^{n_{\beta m}(a) - \beta m}.$

Now for any $\delta' > 0$ and $m$ large enough,

$$\inf_{y \in \partial D(x, r_{m,\beta m})} \mathbf{P}^y(T_z < T_{\partial D(x, r_{m,\beta m-1})})$$

$$\geq \inf_{y \in \partial D(z, r_{m,\beta m-1} + |z-x|)} \mathbf{P}^y(T_z < T_{\partial D(z, r_{m,\beta m-1} - |z-x|)})$$

$$= \frac{\log((r_{m,\beta m-1} - |z-x|)/(r_{m,\beta m-1} + |z-x|)) + O((m \log m)^{-1})}{\log(r_{m,\beta m-1} - |z-x|)} \geq \frac{1 - \delta'}{\gamma m}$$

where we have used Proposition 1.6.7 of [7] in the latter equality. Hence, by (2.15),

$$(2.16) \quad \mathbf{P}(\mathcal{L}_{z,m,\beta}) \leq \left(1 - \frac{1 - \delta'}{\gamma m}\right)^{n_{\beta m}(a) - \beta m} \leq e^{-(1 - 2\delta')3a\beta^2 m (\log m)/\gamma}.$$

Since $r_{m,\beta m+1}^2 \leq e^{2(1+\delta')3\gamma m \log m}$, taking $\delta'$ sufficiently small we get (2.13) from (2.14), (2.16) and (2.5).

Hereafter we write $N \overset{k}{\sim} n_k$ when $|N - n_k(a)| \leq k$ and for any $x \in \mathbb{Z}^2$ and $\rho < R$ let $\mathcal{G}^x(R; \rho)$ denote the $\sigma$-algebra generated by the excursions of the SRW from $\partial D(x, R)$ to $\partial D(x, \rho)$, including the part of the path till first hitting $\partial D(x, R)$. Conditioning on $N_{m,\beta m}^x = \ell$ and on $\mathcal{G}_{\beta m}^x := \mathcal{G}^x(r_{m,\beta m-1}; r_{m,\beta m})$, for each $x \in \mathcal{A}_m$ the event $H_{\beta m}^x$ holds if and only if the SRW visits each site in $D(x, r_{m,\beta m+1})$ during its first $\ell$ excursions from $\partial D(x, r_{m,\beta m})$ to $\partial D(x, r_{m,\beta m-1})$. Since both $r_{m,\beta m+1}/r_{m,\beta m}$ and $r_{m,\beta m}/r_{m,\beta m-1}$ are of $O(m^{-3})$ while $m^{-3}(\log m)n_{\beta m}(a) \to 0$, it follows from Lemma 2.4 of [4] that uniformly with respect to $\ell \overset{\beta m}{\sim} n_{\beta m}$ and $x \in \mathcal{A}_m$

$$(2.17) \quad \mathbf{P}(H_{\beta m}^x \mid \mathcal{G}_{\beta m}^x, N_{m,\beta m}^x = \ell) \geq (1 + o(1_m))(1 - \xi_m)\mathbb{1}_{\{N_{m,\beta m}^x = \ell\}}.$$

Since $\{N_{m,k}^x \overset{k}{\sim} n_k\} \in \mathcal{G}_{\beta m}^x$ for any $3 \leq k \leq \beta m$, we deduce by (2.17) that

$$\mathbf{P}(x \text{ is } m, \beta\text{-successful})$$

$$(2.18) \qquad = \mathbf{P}(N_{m,k}^x \overset{k}{\sim} n_k \ \forall k \in [3, \beta m]; H_{\beta m}^x)$$

$$= \sum_{\ell \overset{\beta m}{\sim} n_{\beta m}} \mathbf{P}(N_{m,k}^x \overset{k}{\sim} n_k \ \forall k \in [3, \beta m-1]; N_{m,\beta m}^x = \ell; H_{\beta m}^x)$$

$$\geq (1 + o(1_m))(1 - \xi_m)\mathbf{P}(N_{m,k}^x \overset{k}{\sim} n_k \ \forall k \in [3, \beta m]).$$

Recall that $\xi_m \to 0$ by the estimate of (2.13), hence

$$\mathbf{P}(x \text{ is } m, \beta\text{-successful}) \geq (1 + o(1_m))\mathbf{P}(N_{m,k}^x \overset{k}{\sim} n_k \ \forall k \in [3, \beta m]),$$

which amounts to (2.12). $\quad\square$

**3. The upper bound for Theorem 1.1.** In this section we establish the upper bound for (2.1):

$$
\limsup_{m \to \infty} \frac{\log \mathcal{R}(r_m)}{\log r_m} \leq \frac{1}{2} \qquad \text{a.s.}
\tag{3.1}
$$

Fix $0 < \gamma < 1$. We begin by describing a two-tiered collection of discs in $D(0, r_m)$. Let $\mathcal{B}_{m,2}$ be a maximal collection of points in $D(0, r_{m-1})$ such that the discs $\{D(x_i, r_{\gamma m+2})\,; x_i \in \mathcal{B}_{m,2}\}$ are disjoint and do not intersect $D(0, r_{\gamma m+2})$ (hereafter we use $\gamma m$ also for the integer part of $\gamma m$). For each $x \in \mathcal{B}_{m,2}$, let $\mathcal{B}_{m,1}(x)$ be a maximal collection of points such that the discs $\{D(y_i, r_{\gamma m})\,; y_i \in \mathcal{B}_{m,1}(x)\}$ are disjoint and contained in $D(x, r_{\gamma m+2})$. Let $\mathcal{B}_{m,1} = \bigcup_{x \in \mathcal{B}_{m,2}} \mathcal{B}_{m,1}(x)$.

For any $y \in \mathcal{B}_{m,1}$, we let $N_m^y$ denote the number of excursions from $\partial D(y, r_{\gamma m-1})$ to $\partial D(y, r_{\gamma m})$ until time $\tau(r_m)$. Recall the notation $n_k(a) = 3ak^2 \log k$ for $a > 0$, $k \geq 3$, and taking $\beta = 1 - \gamma$, consider the events

$$
\Gamma_m(a) = \bigcap_{y \in \mathcal{B}_{m,1}} \{N_m^y \leq n_{\beta m}(a)\},
\tag{3.2}
$$

about which the following lemma is proved at the end of the section.

LEMMA 3.1. *For any $a > 2$ we can find $\zeta = \zeta(a, \gamma) > 0$ such that for all $m$ sufficiently large,*

$$
\mathbf{P}(\Gamma_m(a)) \geq 1 - e^{-\zeta m \log m}.
\tag{3.3}
$$

For any $D(x, r) \subseteq \mathbb{Z}^2$, let $\mathcal{C}(x, r)$ denote the number of steps it takes the SRW to cover $D(x, r)$. Note that if the SRW covers a disc of radius $r_{\gamma m+3}$ with center in $D(0, r_{m-1})$, then it must also cover $D(x, r_{\gamma m+2})$ for some $x \in \mathcal{B}_{m,2}$. Therefore, one can easily check that the upper bound (3.1) follows once we show that for any $1/2 < \gamma < 1$

$$
\sum_{m=5}^{\infty} \mathbf{P}\left( \bigcup_{x \in \mathcal{B}_{m,2}} \{\mathcal{C}(x, r_{\gamma m+2}) \leq \tau(r_m)\} \right) < \infty.
\tag{3.4}
$$

Further, by Lemma 3.1 it suffices to show that for some $a > 2$,

$$
\sum_{m=5}^{\infty} \mathbf{P}\left( \bigcup_{x \in \mathcal{B}_{m,2}} \{\mathcal{C}(x, r_{\gamma m+2}) \leq \tau(r_m)\} \mid \Gamma_m(a) \right) < \infty.
\tag{3.5}
$$

We have

$$
\begin{aligned}
\mathbf{P}\Bigg( &\bigcup_{x \in \mathcal{B}_{m,2}} \{\mathcal{C}(x, r_{\gamma m+2}) \leq \tau(r_m)\} \mid \Gamma_m(a) \Bigg) \\
&\leq \sum_{x \in \mathcal{B}_{m,2}} \mathbf{P}(\mathcal{C}(x, r_{\gamma m+2}) \leq \tau(r_m) \mid \Gamma_m(a)),
\end{aligned}
\tag{3.6}
$$

and since $|\mathcal{B}_{m,2}| \leq Ce^{m^2}$, for $\gamma > 1/2$, choosing $a(\gamma) > 2$ sufficiently close to 2, the summability of (3.5) is a consequence of the next lemma.

LEMMA 3.2.   *If* $0 < \gamma = 1 - \beta < 1$ *and* $a > 0$ *are such that*

$$(3.7) \qquad a\beta^2 < 2\gamma^2,$$

*then for all $m$ large enough,*

$$(3.8) \qquad \sup_{x \in \mathcal{B}_{m,2}} \mathbf{P}(\mathcal{C}(x, r_{\gamma m + 2}) \leq \tau(r_m) \mid \Gamma_m(a)) \leq e^{-m^3}.$$

PROOF OF LEMMA 3.2.   Clearly,

$$(3.9) \qquad \{\mathcal{C}(x, r_{\gamma m + 2}) \leq \tau(r_m)\} \subseteq \bigcap_{y \in \mathcal{B}_{m,1}(x)} \{\mathcal{C}(y, r_{\gamma m - 2}) \leq \tau(r_m)\}.$$

Let $\widetilde{\mathcal{C}}(y, r_{\gamma m - 2}; k)$ denote the event that $D(y, r_{\gamma m - 2})$ is covered in the first $k$ excursions from $\partial D(y, r_{\gamma m - 1})$ to $\partial D(y, r_{\gamma m})$. Note that for any $y \in \mathcal{B}_{m,1}$, the events $\Gamma_m(a)$ and $\{\mathcal{C}(y, r_{\gamma m - 2}) \leq \tau(r_m)\}$ imply that also $\widetilde{\mathcal{C}}(y, r_{\gamma m - 2}; n_{\beta m}(a))$ holds. Our construction of $\mathcal{B}_{m,1}$ guarantees that for each $m$, $k$ and $y', y \in \mathcal{B}_{m,1}$ such that $y' \neq y$, the events $\Gamma_m(a)$ and $\widetilde{\mathcal{C}}(y', r_{\gamma m - 2}; k)$ are in the $\sigma$-algebra $\mathcal{G}^y(r_{\gamma m}; r_{\gamma m - 1})$ generated by the excursions of the SRW from $\partial D(y, r_{\gamma m})$ to $\partial D(y, r_{\gamma m - 1})$. Further, by Lemma 2.4 of [4] we know that uniformly in $y \in \mathcal{B}_{m,1}$

$$\mathbf{P}(\widetilde{\mathcal{C}}(y, r_{\gamma m - 2}; n_{\beta m}(a)) \mid \mathcal{G}^y(r_{\gamma m}; r_{\gamma m - 1})) = (1 + o(1_m))\mathbf{P}(\widetilde{\mathcal{C}}(y, r_{\gamma m - 2}; n_{\beta m}(a))).$$

Consequently, by (3.9),

$$\mathbf{P}(\mathcal{C}(x, r_{\gamma m + 2}) \leq \tau(r_m) \mid \Gamma_m(a)) \leq \mathbf{P}\left(\bigcap_{y \in \mathcal{B}_{m,1}(x)} \widetilde{\mathcal{C}}(y, r_{\gamma m - 2}; n_{\beta m}(a)) \mid \Gamma_m(a)\right)$$

$$= \prod_{y \in \mathcal{B}_{m,1}(x)} (1 + o(1_m))\mathbf{P}(\widetilde{\mathcal{C}}(y, r_{\gamma m - 2}; n_{\beta m}(a))).$$

We will show that if $a$, $\gamma$ and $\beta = 1 - \gamma$ satisfy (3.7), then

$$(3.10) \qquad \sup_{y \in \mathcal{B}_{m,1}} \mathbf{P}(\widetilde{\mathcal{C}}(y, r_{\gamma m - 2}; n_{\beta m}(a))) = o(1_m).$$

Since $|\mathcal{B}_{m,1}(x)| \geq m^4$, this in turn results with the statement (3.8) of the lemma.

To prove (3.10) we fix $\gamma' > \gamma$ and applying (3.19) of [4] with $K = r_{\gamma' m}$, $R = r_{\gamma m}$, $r = r_{\gamma m - 1}$ and $N = n_{\beta m}(a)$, we deduce that for any $\delta > 0$, uniformly in $y \in \mathcal{B}_{m,1}$ with probability $1 - o(1_m)$ it takes the SRW on the two dimensional torus $\mathbb{Z}_K^2$ of side length $K$ less than $T := \frac{2}{\pi}(1 + \delta)K^2 N \log(R/r)$ steps to complete $N$ excursions from $\partial D(y, r_{\gamma m-1})$ to $\partial D(y, r_{\gamma m})$ [since $N \log(R/r)/\log(K/r) \to \infty$ as $m \to \infty$]. It is not hard to verify that (3.7) implies that $T \leq \frac{4}{\pi}(1-\delta)(K \log K)^2$ for $\delta = \delta(a, \gamma) > 0$ sufficiently small and all $m$ large enough. It then follows from (1.2) of [4] that when $\gamma'$ is such that $(\gamma/\gamma')^2 > (1-\delta)$, the probability that the SRW on the torus $\mathbb{Z}_K^2$ covers $D(y, r_{\gamma m-2})$ within that many steps is $o(1_m)$, again uniformly in $y$. Thus, the probability that the SRW on the torus $\mathbb{Z}_K^2$ covers $D(y, r_{\gamma m-2})$ during its first $N$ excursions from $\partial D(y, r_{\gamma m-1})$ to $\partial D(y, r_{\gamma m})$ is $o(1_m)$. We are interested in the same probability, but for the SRW on $\mathbb{Z}^2$. However, note that conditioned on their beginning and end points the $N$ excursions in question are mutually independent and each has the same (conditioned) law for $\mathbb{Z}^2$ and for the torus $\mathbb{Z}_K^2$. Further, from Lemma 2.4 of [4] we know that the probability we are considering is, up to a factor $1 + o(1_m)$, independent of the beginning and end points of these excursions. This completes the proof of (3.10) and hence of the lemma. $\quad\square$

PROOF OF LEMMA 3.1. The proof is similar to that of (2.13). Indeed, fixing $a > 2$ it suffices to show that for some $\zeta = \zeta(a, \gamma) > 0$ and all $m$ large enough

$$(3.11) \qquad \mathbf{P}\left(\bigcup_{y \in \mathcal{B}_{m,1}} \{N_m^y > n_{\beta m}(a)\}\right) \leq e^{-\zeta m \log m}.$$

To see this, note that

$$(3.12) \quad \mathbf{P}\left(\bigcup_{y \in \mathcal{B}_{m,1}} \{N_m^y > n_{\beta m}(a)\}\right) \leq |\mathcal{B}_{m,1}| \sup_{y \in \mathcal{B}_{m,1}} \mathbf{P}(N_m^y > n_{\beta m}(a))$$

and by the strong Markov property of the SRW,

$$\mathbf{P}(N_m^y > n_{\beta m}(a)) \leq \left(\sup_{x \in \partial D(y, r_{\gamma m})} \mathbf{P}^x(T_{\partial D(0, r_m)} > T_{\partial D(y, r_{\gamma m-1})})\right)^{n_{\beta m}(a)}.$$

(3.13)

Now for any $\delta' > 0$ and $m$ large enough,

$$\inf_{x \in \partial D(y, r_{\gamma m})} \mathbf{P}^x(T_{\partial D(0, r_m)} < T_{\partial D(y, r_{\gamma m-1})})$$

$$\geq \inf_{x \in \partial D(y, r_{\gamma m})} \mathbf{P}^x(T_{\partial D(y, r_m + |y|)} < T_{\partial D(y, r_{\gamma m-1})})$$

$$= \frac{\log(r_{\gamma m}/r_{\gamma m-1}) + O(r_{\gamma m-1}^{-1})}{\log((r_m + |y|)/r_{\gamma m-1})} \geq \frac{1 - \delta'}{\beta m}$$

where we have used Exercise 1.6.8 of [7] in the latter equality. Hence, by (3.13),

$$(3.14) \quad \mathbf{P}(N_m^y > n_{\beta m}(a)) \leq \left(1 - \frac{1 - \delta'}{\beta m}\right)^{n_{\beta m}(a)} \leq e^{-a(1-2\delta')3\beta m \log m}.$$

Since $|\mathcal{B}_{m,1}| \le e^{2(1+\delta')3\beta m \log m}$, for $a > 2$ the estimate (3.11) follows from (3.12) and (3.14) upon taking $\delta' > 0$ sufficiently small. $\square$

**4. Proof of Theorem 1.3.** By the same argument as in the proof of the lower bound for (1.3), the lower bound for (1.8) follows once we show that for any $\eta > 0$ there exists $p_\eta > 0$ such that for all sufficiently large $m$,

$$(4.1) \qquad \mathbf{P}\Big(\frac{\log \mathcal{R}_\ell(8r_m)}{\log r_m} \ge \frac{1}{1+\sqrt{\ell}} - 2\eta\Big) \ge p_\eta > 0.$$

To this end, for $x \in \mathbb{Z}^2$, $3 \le k \le m-1$ and $1 \le j \le \ell$ let $N_{m,k}^{x,j}$ denote the number of excursions of the $j$'th SRW from $\partial D(x, r_{m,k})$ to $\partial D(x, r_{m,k-1})$ until its hitting time of $\partial D(x, r_m)$. Fixing $a > 0$ and $0 < \beta = 1 - \gamma < 1$, we set $n_k(a) = 3ak^2 \log k$. We say that $x \in \mathbb{Z}^2$ is $m, \beta, \ell$-*presuccessful* if

$$n_k(a) - k \le N_{m,k}^{x,j} \le n_k(a) + k \qquad \text{for } k = 3, \dots, \beta m, j = 1, \dots, \ell,$$

and say that $x \in \mathbb{Z}^2$ is $m, \beta, \ell$-*successful* if in addition for $j = 1, \dots, \ell$, each point in $D(x, r_{m,\beta m+1})$ is visited during the first $N_{m,\beta m}^{x,j}$ excursions of the $j$'th SRW from $\partial D(x, r_{m,\beta m})$ to $\partial D(x, r_{m,\beta m-1})$.

As before, $\mathcal{A}_m \subseteq \mathbb{Z}^2$ is a maximal collection of points in $[3r_m, 4r_m]^2$ such that the distance between any two points in $\mathcal{A}_m$ is at least $4r_{m,\beta m}$. With $\log r_{m,\beta m+1} \sim \gamma \log r_m$ we establish (4.1) by showing that for some $a(\beta)$ and any $\gamma = 1 - \beta < 1/(1+\sqrt{\ell})$ the probability that there is at least one $m, \beta, \ell$-successful point in $\mathcal{A}_m$ is bounded away from zero as $m \to \infty$. Specifically, as before we show that if $a$ and $\beta = 1 - \gamma$ satisfy (2.5) then uniformly in $\mathcal{A}_m$ each $m, \beta, \ell$-presuccessful point is with high probability also $m, \beta, \ell$-successful. Then we show that $a < 2/\ell$ guarantees that with a probability that is bounded away from zero as $m \to \infty$, there exists at least one $m, \beta, \ell$-presuccessful point in $\mathcal{A}_m$. This establishes (4.1) because for $\gamma = 1 - \beta < 1/(1+\sqrt{\ell})$ we can satisfy (2.5) with some $a < 2/\ell$. Indeed, arguing as in Section 2, for $a < 2/\ell$ the existence of at least one $m, \beta, \ell$-successful point in $\mathcal{A}_m$ is a consequence of our next lemma [simply take $\eta > 0$ so that $2 - \ell(a + 2\eta) > 0$ when adapting (2.11) to the present context], which thus completes the proof of the lower bound for (1.8).

LEMMA 4.1. *If $a$ and $\beta = 1 - \gamma$ satisfy (2.5) then there exists $\delta_m \to 0$ such that*

$$(4.2) \qquad \bar{q}_{m,\ell} := \inf_{x \in \mathcal{A}_m} \mathbf{P}(x \ \text{is} \ m, \beta, \ell\text{-successful}) \ge r_{\beta m}^{-\ell(a+\delta_m)}.$$

*Further, for some $c < \infty$, all $m$ and $x \ne y \in \mathcal{A}_m$*

$$(4.3) \qquad \mathbf{P}(x, y \ \text{are} \ m, \beta, \ell\text{-successful}) \le c\bar{q}_{m,\ell}^2 r_{k(x,y)}^{\ell(a+\delta_{k(x,y)})}$$

[where as before $k(x,y) = \min\{j \ge 1 : D(x, r_{m,j}+1) \cap D(y, r_{m,j}+1) = \varnothing\} \le \beta m$].

PROOF. This is an easy adaptation of the proof of Lemma 2.1. Indeed, as shown there, (2.5) guarantees that uniformly in $\mathcal{A}_m$ each $m, \beta, \ell$-presuccessful point is with high probability also $m, \beta, \ell$-successful [raising the factor $(1 + o(1_m))(1 - \zeta_m)$ in (2.18) to the $\ell$th power], while since the presuccessful condition now involves $\ell$ independent walks, the probabilities in the statement of Lemma 2.1 are now raised to the $\ell$th power. □

We turn next to the upper bound for (1.8), which amounts to showing that

$$(4.4) \qquad \limsup_{m \to \infty} \frac{\log \mathcal{R}_\ell(r_m)}{\log r_m} \leq \frac{1}{1 + \sqrt{\ell}} \qquad \text{a.s.}$$

To this end, fixing $0 < \gamma = 1 - \beta < 1$ we adapt the argument of Section 3 using the same two-tiered collection of discs in $D(0, r_m)$. For any $y \in \mathcal{B}_{m,1}$ here $N_m^{y,j}$ denotes the number of excursions of the $j$th SRW from $\partial D(y, r_{\gamma m-1})$ to $\partial D(y, r_{\gamma m})$ until the time $\tau_j(r_m)$ in which the $j$th SRW first exits $D(0, r_m)$. Fixing $a > 0$ and setting again $n_k(a) = 3ak^2 \log k$ we now consider the events

$$(4.5) \qquad \Gamma_{m,\ell}(a) := \bigcap_{y \in \mathcal{B}_{m,1}} \bigcup_{j=1}^{\ell} \{N_m^{y,j} \leq n_{\beta m}(a)\},$$

about which we show the following.

LEMMA 4.2. For any $a > 2/\ell$ we can find $\zeta = \zeta(a, \gamma) > 0$ such that for all $m$ sufficiently large,

$$(4.6) \qquad \mathbf{P}(\Gamma_{m,\ell}(a)) \geq 1 - e^{-\zeta m \log m}.$$

PROOF. It suffices to show that when $a > 2/\ell$ we can find $\zeta = \zeta(a, \gamma) > 0$ such that for all $m$ sufficiently large,

$$(4.7) \qquad \mathbf{P}\left(\bigcup_{y \in \mathcal{B}_{m,1}} \bigcap_{j=1}^{\ell} \{N_m^{y,j} > n_{\beta m}(a)\}\right) \leq e^{-\zeta m \log m}.$$

Since the $\ell$ walks are independent, we have that

$$(4.8) \quad \mathbf{P}\left(\bigcup_{y \in \mathcal{B}_{m,1}} \bigcap_{j=1}^{\ell} \{N_m^{y,j} > n_{\beta m}(a)\}\right) \leq |\mathcal{B}_{m,1}| \left(\sup_{y \in \mathcal{B}_{m,1}} \mathbf{P}(N_m^y > n_{\beta m}(a))\right)^{\ell}$$

and using the upper bound of (3.14) with $\delta'$ sufficiently small, we verify that (4.7) holds for some $\zeta > 0$ as soon as $\ell a > 2$. □

For $j = 1, \ldots, \ell$, let $\mathcal{C}_j(x, r)$ denote the number of steps required until the $j$th SRW covers $D(x, r)$. Note that if all $\ell$ walks cover some disc of radius $r_{\gamma m+3}$ with center in $D(0, r_{m-1})$, then necessarily for some $x \in \mathcal{B}_{m,2}$ they all cover $D(x, r_{\gamma m+2})$. It is therefore easy to check that the upper bound of (4.4) follows once we show that for each $1/(1 + \sqrt{\ell}) < \gamma < 1$,

$$(4.9) \qquad \sum_{m=5}^{\infty} \mathbf{P} \left( \bigcup_{x \in \mathcal{B}_{m,2}} \bigcap_{j=1}^{\ell} \{ \mathcal{C}_j(x, r_{\gamma m+2}) \le \tau_j(r_m) \} \right) < \infty.$$

Further, considering $a \downarrow 2/\ell$, in view of Lemma 4.2 this is a direct consequence of our next lemma.

LEMMA 4.3. *If $0 < \gamma = 1 - \beta < 1$ and $a > 0$ satisfy (3.7), then for all $m$ large enough,*

$$(4.10) \qquad \sup_{x \in \mathcal{B}_{m,2}} \mathbf{P} \left( \bigcap_{j=1}^{\ell} \{ \mathcal{C}_j(x, r_{\gamma m+2}) \le \tau_j(r_m) \} \mid \Gamma_{m,\ell}(a) \right) \le e^{-m^3}.$$

PROOF. Let $\widetilde{\mathcal{C}}_j(y, r_{\gamma m-2}; k)$ denote the event that the $j$th SRW covers $D(y, r_{\gamma m-2})$ during its first $k$ excursions from $\partial D(y, r_{\gamma m-1})$ to $\partial D(y, r_{\gamma m})$. Note that for any $y \in \mathcal{B}_{m,1}$, the events $\Gamma_{m,\ell}(a)$ and $\{\mathcal{C}_j(y, r_{\gamma m-2}) \le \tau(r_m)\}$ for $j = 1, \ldots, \ell$, imply that at least one of the events $\widetilde{\mathcal{C}}_j(y, r_{\gamma m-2}; n_{\beta m}(a))$ holds as well. Thus, by the independence of the $\ell$ walks, as in the proof of Lemma 3.2 we have that for all $x \in \mathcal{B}_{m,2}$,

$$\mathbf{P} \left( \bigcap_{j=1}^{\ell} \{ \mathcal{C}_j(x, r_{\gamma m+2}) \le \tau_j(r_m) \} \mid \Gamma_{m,\ell}(a) \right)$$

$$\le \mathbf{P} \left( \bigcap_{y \in \mathcal{B}_{m,1}(x)} \bigcup_{j=1}^{\ell} \widetilde{\mathcal{C}}_j(y, r_{\gamma m-2}; n_{\beta m}(a)) \mid \Gamma_{m,\ell}(a) \right)$$

$$\le \prod_{y \in \mathcal{B}_{m,1}(x)} \left[ (1 + o(1_m)) \sum_{j=1}^{\ell} \mathbf{P}(\widetilde{\mathcal{C}}_j(y, r_{\gamma m-2}; n_{\beta m}(a))) \right].$$

Since $a$ and $\gamma = 1 - \beta$ satisfy (3.7), we next apply the bound (3.10) to the preceding inequality, and with $|\mathcal{B}_{m,1}(x)| \ge m^4$, thus establish the bound (4.10) of the lemma. □

**5. The lower bound for Theorem 1.2.** As seen before, it suffices to consider the lower bound for (1.6) and the sequence $r_m$. Further, by the same argument as in the proof of the lower bound for Theorem 1.1, it suffices to prove the analog of (2.3), namely, to show that

$$(5.1) \qquad \liminf_{m \to \infty} \mathbf{P}(\log \mathcal{R}(8r_m; 4\alpha (\log r_m)^2/\pi) \ge \gamma \log r_m) > 0,$$

whenever $\gamma = 1 - \beta < (1 - \sqrt{\alpha})/2$ or equivalently, $0 < \gamma < \beta - \sqrt{\alpha}$. To prove (5.1), rerun the arguments of Section 2 while replacing the event $H_{\beta m}^x$ in the definition (2.4) of an $m, \beta$-successful point $x$ with the event $H_{\beta m}^x(\alpha)$ that the SRW visits each point in $D(x, r_{m,\beta m+1})$ at least $\alpha(\log r_m)^2/\pi$ times during its first $N_{m,\beta m}^x$ excursions from $\partial D(x, r_{m,\beta m})$ to $\partial D(x, r_{m,\beta m-1})$. It is not hard to check that this strategy works as soon as (2.12) applies for this definition of $m, \beta$-successful points and some $a = a(\beta) < 2$, whenever $0 < \gamma < \beta - \sqrt{\alpha}$. As for the latter, let $\mathcal{L}_{z,m,\beta}(\alpha)$ denote the event that $z \in D(x, r_{m,\beta m+1})$ is visited less than $4\alpha(\log r_m)^2/\pi$ times during the first $n_{\beta m}(a) - \beta m$ excursions of the SRW from $\partial D(x, r_{m,\beta m})$ to $\partial D(x, r_{m,\beta m})$. Then, following the proof of Lemma 2.1, we have (5.1) as soon as we show that for $0 < \gamma < \beta - \sqrt{\alpha}$ and $a < 2$ sufficiently close to 2,

$$(5.2) \quad \xi_m(\alpha) = \sup_{x \in \mathcal{A}_m} \mathbf{P}\left( \bigcup_{z \in D(x, r_{m,\beta m+1})} \mathcal{L}_{z,m,\beta}(\alpha) \right) \to 0 \qquad \text{as } m \to \infty$$

[compare with (2.13)]. Turning to the derivation of (5.2), set $R = r_{m,\beta m-1} - r_{m,\beta m+1}$ and $\rho = r_{m,\beta m} + r_{m,\beta m+1}$ and let $\mathcal{L}'_{z,m,\beta}(\alpha)$ denote the event that $z$ is visited less than $4\alpha(\log r_m)^2/\pi$ times during the first $n_{\beta m}(a) - \beta m$ excursions of the SRW from $\partial D(z, \rho)$ to $\partial D(z, R)$. Note that if $z \in D(x, r_{m,\beta m+1})$ then

$$D(x, r_{m,\beta m}) \subseteq D(z, \rho) \subseteq D(z, R) \subseteq D(x, r_{m,\beta m-1}),$$

implying that the SRW makes at least $k$ excursions from $\partial D(z, \rho)$ to $\partial D(z, R)$ during its first $k$ excursions from $\partial D(x, r_{m,\beta m})$ to $\partial D(x, r_{m,\beta m-1})$, so in particular, $\mathcal{L}_{z,m,\beta}(\alpha) \subseteq \mathcal{L}'_{z,m,\beta}(\alpha)$. Consequently, with $R > r_{m,\beta m+1}$, for some constant $c < \infty$

$$
\begin{aligned}
(5.3) \quad \xi_m(\alpha) &\le c r_{m,\beta m+1}^2 \sup_{\substack{z \in D(x, r_{m,\beta m+1}) \\ x \in \mathcal{A}_m}} \mathbf{P}(\mathcal{L}_{z,m,\beta}(\alpha)) \\
&\le c R^2 \sup_{z \notin D(0,\rho)} \mathbf{P}(\mathcal{L}'_{z,m,\beta}(\alpha)).
\end{aligned}
$$

When $0 < \gamma = 1 - \beta < \beta - \sqrt{\alpha}$, if both $\delta > 0$ and $2 - a > 0$ are sufficiently small, then our next lemma shows that $\mathbf{P}(\mathcal{L}'_{z,m,\beta}(\alpha)) \le R^{-2-\eta}$ for some $\eta = \eta(\beta, \alpha, \delta, a) > 0$ and all $z \notin D(0, \rho)$. Combining this bound with (5.3) yields that (5.2) holds, thus completing the proof of the lower bound for Theorem 1.2.

LEMMA 5.1.   *If $a, \delta, \beta > 0$ are such that*

$$(5.4) \qquad\qquad (1-\delta)^2 a \beta^2 > 2\alpha,$$

*then for all $m$ sufficiently large,*

$$(5.5) \qquad \sup_{z \notin D(0,\rho)} \mathbf{P}(\mathcal{L}'_{z,m,\beta}(\alpha)) \leq R^{-((1-\delta)\beta\sqrt{a}-\sqrt{2\alpha})^2/\gamma^2}.$$

PROOF. Let $L(j,z)$ denote the number of visits to $z$ during the $j$th excursion of the SRW from $\partial D(z,\rho)$ to $\partial D(z,R)$. Setting $k = k(a, \beta m) := n_{\beta m}(a) - \beta m$ and $s = s(\alpha, m) := 4\alpha(\log r_m)^2/\pi$, by Chebyshev's inequality and the strong Markov property of the SRW, for any $\lambda > 0$

$$\mathbf{P}(\mathcal{L}'_{z,m,\beta}(\alpha)) = \mathbf{P}\left(\sum_{j=1}^{k} L(j,z) < s\right)$$

$$(5.6)$$

$$\leq e^{\lambda s}\mathbb{E}(e^{-\lambda \sum_{j=1}^{k} L(j,z)}) \leq e^{\lambda s}\left[\sup_{y \in \partial D(z,\rho)} \mathbb{E}^y(e^{-\lambda L(z)})\right]^k,$$

where $L(z)$ denotes the number of visits to $z$ of a SRW that starts at $y \in \partial D(z,\rho)$ and is killed upon reaching $\partial D(z,R)$. Since the preceding bound is independent of $z$, we take hereafter $z = 0$ and let

$$G_R(v,u) := \mathbb{E}^v\left(\sum_{i=0}^{\tau(R)} \mathbb{1}_{\{S_i = u\}}\right)$$

denote the Green function for the SRW on $D(0,R)$. Clearly, $G_R(y,0) = \mathbb{E}^y L(0)$ and conditional on hitting the origin, $L(0)$ is a geometric random variable. Consequently, $\mathbf{P}^y(L(0) = j + 1) = pq^j(1-q)$ for $j = 0, 1, \ldots$ with $p = G_R(y,0)/G_R(0,0)$ and $q = 1 - 1/G_R(0,0)$. Hence for any $\lambda > 0$,

$$(5.7) \qquad \mathbb{E}^y(e^{-\lambda L(0)}) = 1 - \frac{(e^\lambda - 1)G_R(y,0)}{1 + (e^\lambda - 1)G_R(0,0)}.$$

By Proposition 1.6.6 of [7], $G_R(0,0) \sim \frac{2}{\pi}\log R$ when $R \to \infty$, so taking $\lambda = \frac{\pi}{2}\varphi/\log R$ we have that $(e^\lambda - 1)G_R(0,0) \sim \varphi$ as $m \to \infty$ (i.e., $R \to \infty$). Further, by Proposition 1.6.7 of [7] we have that

$$\inf_{y \in \partial D(0,\rho)} G_R(y,0) = \frac{2}{\pi}\log\left(\frac{R}{\rho}\right) + O(\rho^{-1}).$$

Recall that our choices of $R = r_{m,\beta m - 1} - r_{m,\beta m + 1}$ and $\rho = r_{m,\beta m} + r_{m,\beta m + 1}$ are such that $\log(R/\rho)/\log R \sim 1/(\gamma m)$, so we get from (5.7) that for $\lambda = \frac{\pi}{2}\varphi/\log R$, any $\varphi, \delta > 0$ and all $m$ large enough,

$$(5.8) \qquad \sup_{y \in \partial D(0,\rho)} \mathbb{E}^y(e^{-\lambda L(0)}) \leq 1 - \frac{(1-\delta)}{\gamma m}\frac{\varphi}{1+\varphi} \leq e^{-(1-\delta)\varphi/(\gamma m(1+\varphi))}.$$

Hence, by (5.6), for any $\varphi > 0$,

$$(5.9) \qquad \mathbf{P}(\mathcal{L}'_{z,m,\beta}(\alpha)) \leq R^{A\varphi - B\varphi/(1+\varphi)},$$

where by (5.4), for all $m$ large enough,

$$(5.10) \quad A := \frac{\pi}{2} \frac{s(\alpha, m)}{(\log R)^2} = 2\alpha \left( \frac{\log r_m}{\log R} \right)^2 \sim \frac{2\alpha}{\gamma^2},$$

$$(5.11) \quad B := (1 - \delta) \frac{k(a, \beta m)}{\gamma m \log R} = (1 - \delta) \frac{n_{\beta m}(a) - \beta m}{\gamma m \log R} \sim \frac{(1 - \delta) a \beta^2}{\gamma^2}$$

are such that $B > A > 0$, in which case a straightforward computation shows that

$$(5.12) \quad \inf_{\varphi > 0} \left( A\varphi - B \frac{\varphi}{1 + \varphi} \right) = -(\sqrt{B} - \sqrt{A})^2.$$

Combining (5.9)–(5.12) we get that (5.5) holds for all $m$ large enough. $\quad\square$

**6. The upper bound for Theorem 1.2.** As explained before, it suffices to prove the upper bound in (1.6) for the sequence $r_m$. That is, fixing $1/2 > \gamma > (1 - \sqrt{\alpha})/2$, to show that

$$(6.1) \quad \limsup_{m \to \infty} \frac{\log \mathcal{R}(r_m; 4\alpha(\log r_m)^2/\pi)}{\log r_m} < \gamma \qquad \text{a.s.}$$

More precisely, adapting the proof of the upper bound for Theorem 1.1, we show that for such $\gamma$ any disc of radius $r_{\gamma m+3}$ with center in $D(0, r_{m-1})$ contains at time $\tau(r_m)$ sites which the SRW visited less than $4\alpha(\log r_m)^2/\pi$ times. To this end, let $\widetilde{\mathcal{C}}_\alpha(y, r_{\gamma m-2}; k)$ denote the event that every point in $D(y, r_{\gamma m-2})$ is visited at least $4\alpha(\log r_m)^2/\pi$ times during the first $k$ excursions from $\partial D(y, r_{\gamma m-1})$ to $\partial D(y, r_{\gamma m})$. Using the two-tiered collection of discs as in Section 3, upon applying Lemma 3.1 and adapting to the present context the reasoning which precedes (3.10), we find that it suffices to prove the following lemma.

LEMMA 6.1. *If $\beta = 1 - \gamma$, $1/2 > \gamma > \beta - \sqrt{\alpha}$ and $a > 2$ is sufficiently close to 2 for*

$$(6.2) \quad \sqrt{a}\beta - \sqrt{2}\sqrt{\alpha} < \sqrt{2}\gamma,$$

*then*

$$(6.3) \quad \sup_{y \in \mathcal{B}_{m,1}} \mathbf{P}(\widetilde{\mathcal{C}}_\alpha(y, r_{\gamma m-2}; n_{\beta m}(a))) = o(1_m).$$

PROOF. In view of (6.2) we fix $0 < \eta < \gamma$ and $1 < h < 2$ such that

$$\sqrt{a}\beta - \sqrt{2}\sqrt{\alpha} < \sqrt{h}(\gamma - 2\eta).$$

Setting $A = \sqrt{a/h}\beta - \gamma > 0$ let $\widehat{n}_k(h) = 3h(k + Am)^2 \log m$ for $k = 1, \ldots, \gamma m$, noting that

$$(6.4) \quad \widehat{n}_{\gamma m}(h) = 3a\beta^2 m^2 \log m \geq n_{\beta m}(a),$$

and further, for some $a' < 2\alpha$ and all $m$ large enough,

$$(6.5) \qquad\qquad \widehat{n}_{\eta m}(h) + \eta m \le 3a' m^2 \log m$$

[e.g., $a' = h(2\eta + A)^2$ will do].

Next, let $N^z_{\gamma m,k}$ for $k = 1, \ldots, \gamma m - 2$, denote the number of excursions of the SRW from $\partial D(z, r_{k-1})$ to $\partial D(z, r_k)$ during its first $\widehat{n}_{\gamma m}(h)$ excursions from $\partial D(z, \widehat{\rho})$ to $\partial D(z, \widehat{R})$, where $\widehat{R} = r_{\gamma m} + r_{\gamma m - 2}$ and $\widehat{\rho} = r_{\gamma m - 1} - r_{\gamma m - 2}$. We say that $z \notin D(0, \widehat{\rho})$ is $m, \gamma$-presluggish if

$$\widehat{n}_k(h) - k \le N^z_{\gamma m,k} \le \widehat{n}_k(h) + k, \qquad \text{for } k = \eta m, \ldots, \gamma m - b,$$

for some fixed $b \ge 4$ to be determined in the sequel. An $m, \gamma$-presluggish point $z$ is called $m, \gamma$-sluggish if during the first $3a' m^2 \log m$ excursions of the SRW from $\partial D(z, r_{\eta m - 1})$ to $\partial D(z, r_{\eta m})$, it visits $z$ less than $4\alpha (\log r_m)^2 / \pi$ times, an event we denote hereafter by $\widehat{\mathcal{L}}_{z,m,\eta}(\alpha)$.

Note that if $z \in D(y, r_{\gamma m - 2})$ then

$$D(z, \widehat{\rho}) \subseteq D(y, r_{\gamma m - 1}) \subseteq D(y, r_{\gamma m}) \subseteq D(z, \widehat{R}),$$

so prior to completing its first $\widehat{n}_{\gamma m}(h)$ excursions from $\partial D(z, \widehat{\rho})$ to $\partial D(z, \widehat{R})$, the SRW completes that many excursions from $\partial D(y, r_{\gamma m - 1})$ to $\partial D(y, r_{\gamma m})$. Consequently, in view of (6.4) and (6.5), any $m, \gamma$-sluggish point in $D(y, r_{\gamma m - 2})$ is visited by the SRW less than $4\alpha (\log r_m)^2 / \pi$ times during its first $n_{\beta m}(a)$ excursions from $\partial D(y, r_{\gamma m - 1})$ to $\partial D(y, r_{\gamma m})$.

We thus complete the proof of Lemma 6.1 by showing that uniformly in $y \in \mathcal{B}_{m,1}$, with probability $1 - o(1_m)$ there exists an $m, \gamma$-sluggish point in any maximal set $\mathcal{Z}_{\eta m}(y)$ of $4 r_{\eta m}$-separated points in $D(y, r_{\gamma m - 2})$. The key for this is our next lemma (whose proof is deferred to the end of the section).

Lemma 6.2. *There exists $\delta_m \to 0$ such that*

$$(6.6) \qquad \widehat{q}_m := \inf_{z \notin D(0,\widehat{\rho})} \mathbf{P}(z \text{ is } m, \gamma\text{-sluggish}) \ge r_m^{-(\gamma-\eta)h - \delta_m}$$

*and*

$$(6.7) \qquad \sup_{z \notin D(0,\widehat{\rho})} \mathbf{P}(z \text{ is } m, \gamma\text{-sluggish}) = (1 + o(1_m))\widehat{q}_m.$$

*Further, let $k(z, z') = \max\{j : D(z, r_j + 1) \cap D(z', r_j + 1) = \varnothing\}$. Then, for any $\varepsilon > 0$ there exist $C, \kappa < \infty$ which are both independent of $b$, such that for all $m$, and $z, z' \notin D(0, \widehat{\rho})$ with $\eta m \le k(z, z') \le \gamma m - b$*

$$(6.8) \qquad \mathbf{P}(z, z' \text{ are } m, \gamma\text{-sluggish}) \le \widehat{q}_m^2 m^\kappa C^{\gamma m - k(z,z')} \left( \frac{r_{\gamma m - b}}{r_{k(z,z')}} \right)^{h + \varepsilon}.$$

*Furthermore, if $\gamma m - b < k(z, z')$ and $|z - z'| \le 2 r_{\gamma m - 2}$ then*

$$(6.9) \qquad \mathbf{P}(z, z' \text{ are } m, \gamma\text{-sluggish}) \le \widehat{q}_m^2 (1 + o(1_m)).$$

Since there are $r_m^{2(\gamma-\eta)+o(1_m)}$ sites in $\mathcal{Z}_{\eta m}(y)$ and $h < 2$ it follows from (6.6) that the mean number of $m, \gamma$-sluggish points in $\mathcal{Z}_{\eta m}(y)$ diverges as $m \to \infty$. In view of (6.7) and Chebyshev's inequality, we complete the proof of Lemma 6.1 by showing that the second moment of this random variable is $(1 + o(1_m))|\mathcal{Z}_{\eta m}(y)|^2 \hat{q}_m^2$, hence in any disc $D(y, r_{\gamma m-2})$ the probability of finding at least one $m, \gamma$-sluggish point is $1 - o(1_m)$. To this end, by the bound (6.9) it suffices to consider the contribution to the second moment by $z, z' \in \mathcal{Z}_{\eta m}(y)$ with $k = k(z, z') \leq \gamma m - b$. There are at most

$$c_2 |\mathcal{Z}_{\eta m}(y)| \left( \frac{r_{k+1}}{r_{\eta m}} \right)^2 \leq c_3 |\mathcal{Z}_{\eta m}(y)|^2 \left( \frac{r_k}{r_{\gamma m-3}} \right)^2$$

such pairs per given $k$. Hence, by (6.8), for $0 < \varepsilon < 2 - h$ and $b > 3 + \kappa/6$, the contribution of all such pairs to the second moment is at most $|\mathcal{Z}_{\eta m}(y)|^2 \hat{q}_m^2$ times

$$\left( \frac{r_{\gamma m-b}}{r_{\gamma m-3}} \right)^2 c_3 m^\kappa \sum_{k=\eta m}^{\gamma m-b} C^{\gamma m-k} \left( \frac{r_{\gamma m-b}}{r_k} \right)^{h+\varepsilon-2}$$

$$\leq c_4 m^{\kappa-6(b-3)} \sum_{j=0}^{\infty} C^j m^{-3j(2-h-\varepsilon)} = o(1_m),$$

as required for completing the proof. $\square$

PROOF OF LEMMA 6.2. We first show that an $m, \gamma$-presluggish point is with very high probability also $m, \gamma$-sluggish. More precisely, adapting the proof of Lemma 5.1, we shall show that for any $a' < 2\alpha$ and $\eta > 0$ there exists $\varepsilon = \varepsilon(\alpha, a', \eta) > 0$ such that for all large $m$,

$$(6.10) \qquad \inf_{z \notin D(0,\rho)} \mathbf{P}(\widehat{\mathcal{L}}_{z,m,\eta}(\alpha)) \geq 1 - R^{-\varepsilon},$$

where now $\rho := r_{\eta m-1}$ and $R := r_{\eta m}$. Indeed, with $L(j, z)$ denoting the number of visits to $z$ during the $j$th excursion of the SRW from $\partial D(z, \rho)$ to $\partial D(z, R)$, and $L(0)$ denoting the number of visits to 0 of a SRW that starts at $y \in D(0, \rho)$ and is killed upon reaching $\partial D(0, R)$, taking now $s = s(\alpha, m) := 4\alpha(\log r_m)^2/\pi$, and $k = k(a', m) := 3a'm^2 \log m$, we have by Chebyshev's inequality and the strong Markov property of the SRW, that for any $\lambda > 0$

$$1 - \mathbf{P}(\widehat{\mathcal{L}}_{z,m,\eta}(\alpha)) = \mathbf{P}\left( \sum_{j=1}^{k} L(j, z) \geq s \right) \leq e^{-\lambda s} \left[ \sup_{y \in \partial D(0,\rho)} \mathbb{E}^y(e^{\lambda L(0)}) \right]^k.$$

(6.11)

Here $\log(R/\rho)/\log R \sim 1/(\eta m)$, so for $\lambda = \frac{\pi}{2}\varphi/\log R$ the computation leading to (5.8), yields now that for any $1 > \varphi, \delta > 0$ and all $m$ large enough,

$$(6.12) \qquad \sup_{y \in \partial D(0,\rho)} \mathbb{E}^y(e^{\lambda L(0)}) \leq 1 + \frac{(1+\delta)}{\eta m} \frac{\varphi}{1-\varphi} \leq e^{(1+\delta)\varphi/(\eta m(1-\varphi))}.$$

In view of (6.11), taking $\delta > 0$ small enough so $(1+\delta)a' < 2\alpha$, for $\varphi = 1 - \sqrt{B/A} > 0$ we find that

$$(6.13) \qquad 1 - \mathbf{P}(\widehat{\mathcal{L}}_{z,m,\eta}(\alpha)) \leq R^{-A\varphi + B\varphi/(1-\varphi)} = R^{-(\sqrt{A}-\sqrt{B})^2},$$

where we get (6.10) upon checking that

$$(6.14) \qquad A := \frac{\pi}{2}\frac{s(\alpha,m)}{(\log R)^2} \sim \frac{2\alpha}{\eta^2} \quad \text{and} \quad B := (1+\delta)\frac{k(a',m)}{\eta m \log R} \sim \frac{(1+\delta)a'}{\eta^2},$$

so $a'$ and $\delta$ are such that $A > B > 0$ for all $m$ large enough, as needed for (6.13) and (6.10) to hold.

Following the same outline as of the proof of Lemma 2.1, we employ hereafter arguments that are very similar to those in the proof of Lemma 10.1 of [4]. That is, first note that the probability that $z$ is $m, \gamma$-sluggish depends on $z \notin D(0, \widehat{\rho})$ only via the distribution of the SRW upon first hitting $\partial D(z, \widehat{\rho})$. Since the definition of such points involves only $O(m^2 \log m)$ excursions of the walk, whereas $\widehat{R}/\widehat{\rho} = O(m^3)$ and $\widehat{\rho}/r_{\gamma m-b} \geq O(m^3)$, an application of Lemma 2.4 of [4] shows that the dependence of this probability on $z$ is negligible, as stated in (6.7). Similarly, by (6.10) and the fact that $\widehat{\mathcal{L}}_{z,m,\eta}(\alpha)$ is in the $\sigma$-algebra of all excursions from $\partial D(z, r_{\eta m-2})$ to $\partial D(z, r_{\eta m-1})$ completed by the walk during its first $3a'm^2 \log m$ excursions from $\partial D(z, r_{\eta m-1})$ to $\partial D(z, r_{\eta m})$, yet another application of Lemma 2.4 of [4] shows that

$$\widehat{q}_m = (1 + o(1_m)) \inf_{z \notin D(0,\widehat{\rho})} \mathbf{P}(N^z_{\gamma m,k} \overset{k}{\sim} \widehat{n}_k, k = \eta m, \ldots, \gamma m - b)$$

(compare with the derivation leading to (10.10) of [4]). Due to the dependence of the relevant excursions on their terminal points, $\{N^z_{\gamma m,k}\}$ is not a Markov chain. Nevertheless, applying (5.9) of [4], we find that

$$(6.15) \qquad \begin{aligned} \widehat{q}_m = (1 + o(1_m)) &\sum_{\ell_k \overset{k}{\sim} \widehat{n}_k} \mathbf{P}(N^z_{\gamma m,\gamma m-b} = \ell_{\gamma m-b}) \\ &\times 0 \prod_{k=\eta m}^{\gamma m-b-1} \binom{\ell_{k+1} + \ell_k - 1}{\ell_k} p_k^{\ell_k}(1-p_k)^{\ell_{k+1}}, \end{aligned}$$

for $p_k = \log(k+1)/(\log k + \log(k+1))$. Further, it is not hard to check that for some $c_0 < \infty$ and all $\eta m \leq k \leq \gamma m$, if $\ell_k \overset{2k}{\sim} \widehat{n}_k(h)$ then

$$\left| \frac{\ell_k}{\ell_{k+1}} - 1 + \frac{2}{k+Am} \right| \leq \frac{c_0}{m \log m},$$

and hence for some $c_1 < \infty$ and any such $\ell_k$,

$$(6.16) \qquad \frac{k^{-3h-1}}{c_1\sqrt{\log k}} \leq \binom{\ell_{k+1} + \ell_k - 1}{\ell_k} p_k^{\ell_k}(1-p_k)^{\ell_{k+1}} \leq \frac{c_1 k^{-3h-1}}{\sqrt{\log k}}$$

(cf. (10.11) of [4] or Lemma 7.2 of [2]). A similar polynomial bound applies for $\mathbf{P}(N^z_{\gamma m, k} \overset{k}{\sim} \widehat{n}_k(h))$, for example, when $k = \gamma m - b$, so putting (6.15) and (6.16) together we arrive at the bound (6.6).

Since an $m, \gamma$-sluggish point is also $m, \gamma$-presluggish, it suffices to prove the upper bounds of (6.8) and (6.9) for $\mathbf{P}(z, z'$ are $m, \gamma$-presluggish). To this end, with $b \geq 4$, if $|z - z'| \leq 2r_{\gamma m - 2}$ then $D(z, r_{\gamma m - b + 1}) \subseteq D(z', \widehat{\rho})$. Thus, when $\gamma m - b < k(z, z')$ it is easy to verify that the event $\{z'$ is $m, \gamma$-presluggish$\}$ is in the $\sigma$-algebra $\mathcal{G}^z(r_{\gamma m - b + 1}; r_{\gamma m - b})$. Further, as usual, conditioned on $N^z_{\gamma m, \gamma m - b + 1} = \ell$ the event $\{z$ is $m, \gamma$-presluggish$\}$ is in the $\sigma$-algebra of all excursions from $\partial D(z, r_{\gamma m - b - 1})$ to $\partial D(z, r_{\gamma m - b})$ completed by the walk during its first $\ell$ excursions from $\partial D(z, r_{\gamma m - b})$ to $\partial D(z, r_{\gamma m - b + 1})$. Thus, if $\ell$ of the preceding is not too large, then the dependence of $\{z$ is $m, \gamma$-presluggish$\}$ on $\mathcal{G}^z(r_{\gamma m - b + 1}; r_{\gamma m - b})$ is negligible. More precisely, it is not hard to verify that for large enough $m$,

$$\mathbf{P}(N^z_{\gamma m, \gamma m - b + 1} \geq m^2 (\log m)^2) \leq e^{-m^2 \log m} = o(1_m) \widehat{q}^2_m,$$

and we get (6.9) by an application of Lemma 2.4 of [4]. Finally, the proof of (6.8) for presluggish points is a simple adaptation of the arguments used when proving (10.5) of [4]. $\quad\square$

## 7. Proof of Theorem 1.4.
Recall that $V(n)$ is the number of steps after step $n$ until the SRW $(S_i, i \geq 0)$ in $\mathbb{Z}^2$ visits a previously unvisited site. We first prove the upper bound, that is, fixing $1/20 > \varepsilon > 0$, we show that

$$(7.1) \qquad \limsup_{n \to \infty} \frac{\log V(n)}{\log n} \leq \frac{1}{2} + 10\varepsilon \qquad \text{a.s.}$$

To this end, considering the events $\mathcal{J}_n = \{V(n) > n^{1/2 + 10\varepsilon}\}$ and $\mathcal{K}_n = \bigcap_{m > n} \{\widetilde{\mathcal{R}}(m) < m^{1/4 + \varepsilon}\}$, we shall show that

$$(7.2) \qquad \sum_n \mathbf{P}(\mathcal{J}_n \cap \mathcal{K}_n) < \infty.$$

Then, by the Borel–Cantelli lemma, almost surely, $\mathcal{J}_n \cap \mathcal{K}_n$ occurs for only finitely many values of $n$. From Theorem 1.1 we know that almost surely $\mathcal{K}_n$ occurs for all $n$ large enough, thus implying that $\mathcal{J}_n$ occurs for only finitely many values of $n$, and (7.1) ensues.

Turning to prove (7.2), take $\rho = \rho(n) = n^{1/4 + 2\varepsilon}$ and $R = R(n) = \rho^{1 + \varepsilon}$, and let $H(m)$ denote the event that there exists a site $x \in D(S_m, \rho)$ which is not visited by the SRW up to time $m + R^{2 + \varepsilon}$. With $\mathcal{F}_n = \sigma(S_k, k \leq n)$, considering a uniformly chosen site among those in $D(S_m, \rho)$ that are not visited by the SRW up to time $m$, we have by its Markov property that

$$(7.3) \qquad \mathbf{P}(H(m) | \mathcal{F}_m) \leq 1 - \inf_{y \in D(0, \rho)} \mathbf{P}^y(T_0 < R^{2 + \varepsilon}).$$

By Propositions 1.6.6 and 1.6.7 of [7], for $0 < \varepsilon < 1$ and $\rho = \rho(n)$ large enough,

$$(7.4) \qquad \inf_{y \in D(0,\rho)} \mathbf{P}^y(T_0 < \tau(R)) = \inf_{y \in D(0,\rho)} \frac{G_R(y,0)}{G_R(0,0)} \geq \frac{\varepsilon}{2}.$$

Further,

$$\mathbf{P}^y(T_0 < R^{2+\varepsilon}) \geq \mathbf{P}^y(T_0 < \tau(R)) - \mathbf{P}^y(\tau(R) > R^{2+\varepsilon})$$

and $\mathbf{P}^y(\tau(R) > R^{2+\varepsilon}) \leq R^{-\varepsilon/2}$ for all $R$ large enough and $y \in D(0, R)$ (e.g., see inequality (1.2.1) of [7]). Thus, by (7.3) and (7.4) we deduce that $\mathbf{P}(H(m)|\mathcal{F}_m) \leq 1 - \varepsilon/3$ for all $n$ large enough. Now, let $m(i) = n + 1 + (i - 1)n^{1/2+8\varepsilon}$ for $i = 1, \ldots, n^{2\varepsilon}$, and take $n$ large enough for $\rho \geq (n + n^{1/2+10\varepsilon})^{1/4+\varepsilon}$. Since $m(i) + R^{2+\varepsilon} \leq m(i+1)$, it follows that $H(m(i)) \in \mathcal{F}_{m(i+1)}$, and further

$$\mathcal{J}_n \cap \mathcal{K}_n \subseteq \bigcap_{i=1}^{n^{2\varepsilon}} \{ \widetilde{\mathcal{R}}(m(i)) < \rho, V(m(i)) > R^{2+\varepsilon} \} \subseteq \bigcap_{i=1}^{n^{2\varepsilon}} H(m(i)).$$

Consequently, the bound $\mathbf{P}(H(m)|\mathcal{F}_m) \leq 1 - \varepsilon/3$ implies that

$$\mathbf{P}(\mathcal{J}_n \cap \mathcal{K}_n) \leq \mathbf{P}\left( \bigcap_{i \leq n^{2\varepsilon}} H(m(i)) \right) \leq (1 - \varepsilon/3)^{n^{2\varepsilon}},$$

for all $n$ large enough, which results with (7.2).

Fixing $0 < \varepsilon < 1/20$ we conclude the proof by establishing the lower bound

$$(7.5) \qquad \limsup_{n \to \infty} \frac{\log V(n)}{\log n} \geq \frac{1}{2} - 10\varepsilon \qquad \text{a.s.}$$

To this end, consider the stopping times

$$(7.6) \qquad \tau_k = \inf\{n \geq k : S_n \in D(x, n^{1/4-\varepsilon}) \subseteq (S_i, i \leq n) \text{ for some } x \in \mathbb{Z}^2\},$$

for the filtration $\mathcal{F}_n$. That is, $\tau_k$ is the first time $n \geq k$ for which the SRW is in a disc of radius $n^{1/4-\varepsilon}$ having no previously unvisited sites. By Theorem 1.1, almost surely $\{\mathcal{R}(n) > \mathcal{R}(n-1) > n^{1/4-\varepsilon}\}$ for infinitely many $n$ values, each of which satisfies the conditions of (7.6). Consequently, almost surely $\tau_k < \infty$ for all $k$. Completely ordering $\mathbb{Z}^2$ in agreement with the Euclidean distance from the origin, let $X_k$ denote the site closest to $S_{\tau_k}$ among those $x \in \mathbb{Z}^2$ such that every site in $D(x, \tau_k^{1/4-\varepsilon})$ is visited by the SRW by time $\tau_k$. Then, $S_{\tau_k} \in D(X_k, \tau_k^{1/4-\varepsilon})$ and $X_k$ is measurable with respect to $\mathcal{F}_{\tau_k}$.

We next show that the events

$$(7.7) \qquad \mathcal{M}_k = \{V(n) \geq n^{1/2-6\varepsilon} \text{ for some } \tau_k \leq n \leq 2\tau_k\},$$

are such that for some finite $k_0$,

$$(7.8) \qquad \mathbf{P}(\mathcal{M}_k) \geq \tfrac{1}{6} \qquad \forall k \geq k_0.$$

To this end, let $\theta(i)$ denote the shift of the SRW by $i$, that is, considering $\{S_{n+i}\}$ instead of $\{S_n\}$, and let $T^{(k)}(A)$ denote the first hitting time of a set $A \in \mathbb{Z}^2$ by the shifted random walk $S_n^{(k)} := S_n \circ \theta(\tau_k)$. Consider the stopping times $\sigma_k := T^{(k)}(D(X_k, \tau_k^{1/4-2\varepsilon}))$ with respect to the canonical filtration of the shifted walk $(S_i^{(k)}, i \geq 0)$. Note that $S_0^{(k)} = S_{\tau_k} \in D(X_k, \rho_k^{1-4\varepsilon})$ for $\rho_k = \tau_k^{1/4}$, whereas $\sigma_k = \inf\{i \geq 0 : S_i^{(k)} \in D(X_k, \rho_k^{1-8\varepsilon})\}$ and $\rho_k \geq k^{1/4}$ by the definition of $\tau_k$. Therefore, by the strong Markov property of the SRW at the stopping time $\tau_k$ and considering the worst possible choice of $\rho_k$, $X_k$ and $S_{\tau_k}$, we have that for all $k$ sufficiently large,

$$(7.9) \quad \mathbf{P}(\sigma_k \leq \tau_k | \mathcal{F}_{\tau_k}) \geq \inf_{R \geq k^{1/4}} \inf_x \inf_{y \in D(x, R^{1-4\varepsilon})} \mathbf{P}^y(T_{D(x, R^{1-8\varepsilon})} \leq R^4) \geq \tfrac{1}{3}$$

(using Exercise 1.6.8 of [7] in the rightmost inequality). Similarly, for all $k$ sufficiently large, the events

$$\mathcal{V}_k := \{T^{(k)}(\partial D(S_0^{(k)}, \tau_k^{1/4-2\varepsilon})) \geq \tau_k^{1/2-5\varepsilon}\},$$

are such that

$$\mathbf{P}(\mathcal{V}_k \circ \theta(\sigma_k) | \mathcal{F}_{\tau_k}) \geq \inf_{R \geq k^{1/4}} \mathbf{P}(\tau(R^{1-8\varepsilon}) \geq R^{2(1-8\varepsilon)-4\varepsilon}) \geq \tfrac{5}{6}.$$

The lower bound (7.8) then follows from the inclusion

$$(7.10) \quad \{\sigma_k \leq \tau_k\} \cap \{\mathcal{V}_k \circ \theta(\sigma_k)\} \subseteq \mathcal{M}_k.$$

To see this inclusion, note that $S_{\tau_k+\sigma_k} \in D(X_k, \tau_k^{1/4-2\varepsilon})$ and the event $\mathcal{V}_k \circ \theta(\sigma_k)$ guarantees that it takes $S_{\sigma_k+i}^{(k)} = S_{\sigma_k+\tau_k+i}$ at least $\tau_k^{1/2-5\varepsilon}$ steps to travel a distance of $\tau_k^{1/4-2\varepsilon}$ from its position at $i=0$, a fortiori before exiting the disc $D(X_k, \tau_k^{1/4-\varepsilon})$, all the sites of which have been previously visited by the SRW. Consequently, if also $\sigma_k \leq \tau_k$, then

$$V(\tau_k + \sigma_k) \geq \tau_k^{1/2-5\varepsilon} \geq (\tau_k + \sigma_k)^{1/2-6\varepsilon},$$

hence $\mathcal{M}_k$ holds as well.

Since $\tau_k$ are a.s. finite we can find a deterministic function $\psi(k)$ such that $\mathbf{P}(\tau_k > \psi(k)) \leq 1/18$ for all $k$. Then, by (7.8) the events $\mathcal{I}_k := \mathcal{M}_k \cap \{k \leq \tau_k \leq \psi(k)\}$ are such that $\mathbf{P}(\mathcal{I}_k) \geq 1/9$ for all $k \geq k_0$. With $\mathbf{P}^x(\mathcal{I}_k)$ independent of $x$, we see by the Markov property of the SRW that for any $m$,

$$(7.11) \quad \mathbf{P}(\mathcal{I}_k \circ \theta(m) \mid \mathcal{F}_m) = \mathbf{P}^{S_m}(\mathcal{I}_k) = \mathbf{P}(\mathcal{I}_k) \geq \tfrac{1}{9} \qquad \text{a.s. } \forall k \geq k_0.$$

Define inductively the nonrandom $t_1 = k_0$ and $t_j = t_{j-1} + 3\psi(t_{j-1})$ for $j \geq 2$. Then, by (7.11)

$$(7.12) \quad \sum_{j=2}^{\infty} \mathbf{P}(\mathcal{I}_{t_j} \circ \theta(t_j)) = \infty.$$

With $\mathcal{I}_k \in \mathcal{F}_{3\psi(k)}$, it follows that $\mathcal{I}_{t_j} \circ \theta(t_j) \in \mathcal{F}_{t_j + 3\psi(t_j)} = \mathcal{F}_{t_{j+1}}$. Consequently, by the Markov property and the fact, mentioned above, that $\mathbf{P}^x(\mathcal{I}_k)$ is independent of $x$, the events $\{\mathcal{I}_{t_j} \circ \theta(t_j); j \geq 2\}$ are mutually independent. Therefore, by (7.12) and the second Borel–Cantelli lemma, with probability one, infinitely many of them occur. It follows from (7.7) that $\mathcal{I}_k \circ \theta(k)$ readily implies that $V(n+k) \geq n^{1/2 - 6\varepsilon}$ for some $n \geq k$. Thus, with $t_j \uparrow \infty$, clearly (7.5) follows.

**Acknowledgments.** We are grateful to Zhan Shi and Ofer Zeitouni for useful discussions, and to Olivier Daviaud, J. Ben Hough and the associate editor for suggestions that led to improvement of this manuscript.

A. Dembo
Departments of Mathematics
  and of Statistics
Stanford University
Stanford, California 94305
USA
E-mail: amir@math.stanford.edu

Y. Peres
Departments of Statistics
  and of Mathematics
University of California at Berkeley
Berkeley, California 94720
USA
E-mail: peres@stat.berkeley.edu

J. Rosen
Department of Mathematics
College of Staten Island, CUNY
Staten Island, New York 10314
USA
E-mail: jrosen3@earthlink.net